\documentclass[12pt]{amsart}
\usepackage{amscd,amssymb,amsmath}
\usepackage[arrow,matrix]{xy}
\usepackage{graphicx}
\usepackage{cite}
\usepackage{ifpdf}

\newtheorem{thm}[subsection]{Theorem}

\newtheorem{pr}[subsection]{Proposition}
\newtheorem{cor}[subsection]{Corollary}
\newtheorem{rem}[subsection]{Remark}
\newtheorem{defn}[subsection]{Definition}
\newtheorem{ex}[subsection]{Example}

\numberwithin{equation}{section} 

\newcommand{\bea}{\begin{eqnarray}}
\newcommand{\eea}{\end{eqnarray}}

\begin{document}
\title [ Some properties of evolution algebras]
{ Some properties of evolution algebras}

\author {L. M. Camacho, J. R. G\'omez, B. A. Omirov, R. M. Turdibaev}

\address{L.\ M.\ Camacho and J. R. G\'omez\\ Dpto. Mathem\'atica Aplicada I, Universidad
de Sevilla, Avda. Reina Mercedes, s/n. 41012, Sevilla, Spain}
 \email {lcamacho@us.es, jrgomez@us.es}
 \address{B.\ A.\ Omirov\\ Institute of mathematics and information
 technologies, Do'rmon Yo'li str.,29, 100125,
Tashkent, Uzbekistan.} \email {omirovb@mail.ru}
 \address{R.\ M.\ Turdibaev\\ Department of Mathematics, National University of Uzbekistan, Vuzgorogok, 27, 100174,
Tashkent, Uzbekistan.} \email {rustamtm@yahoo.com}
 \maketitle
\begin{abstract}
The paper is devoted to the study of finite dimensional complex
evolution algebras. The class of evolution algebras isomorphic to
evolution algebras with Jordan form matrices is described. For
finite dimensional complex evolution algebras the criteria of
nilpotency is established in terms of the properties of
corresponding matrices. Moreover, it is proved that for nilpotent
$n-$dimensional complex evolution algebras the possible maximal
nilpotency index is $1+2^{n-1}.$ The criteria of planarity for
finite graphs is formulated by means of evolution algebras defined
by graphs.

{\it AMS classifications (2010):} 05C25; 17A36; 17D92; 17D99\\[2mm]

\end{abstract}
\maketitle
\section{Introduction}

In 20s and 30s of the last century the new object was introduced to mathematics, which was the product of interactions between Mendelian genetics and mathematics. Mendel established the basic laws for inheritance, which are summarized as Mendel's Law of Segregation and Mendel's Law of Independent Assortment. This laws were mathematically formulated by Serebrowsky \cite{SER}, who was also the first to give an algebraic interpretation of the $"\times"$ sign, which indicated sexual reproduction. Later Glivenkov \cite{GLI} used the notion of Mendelian algebras in his work. Also Kostitzin \cite{KOS} independently introduced a "symbolic multiplication" to express Mendel's laws. In his several papers Etherington \cite{ETH1}-\cite{ETH3} introduced the formal language of abstract algebra to the study of the genetics. These algebras, in general, are non-associative. 

However, in the beginning of the XX century in genetics there were discovered several examples of inheritances, where traits do not segregate in accordance with Mendel's laws. In the present day, non-Mendelian genetics is a basic language of molecular genetics. Non-Mendelian inheritance plays an important role in several disease processes. Naturally, the question arises: what non-Mendelian genetics offers to mathematics? The evolution algebras, introduced in \cite{TIAN1} serves as the answer to this question.

The concept of evolution algebras lies between algebras and dynamical systems. Algebraically, evolution algebras are non-associative Banach algebra; dynamically, they represent discrete dynamical systems. Evolution algebras have many connections with various branches of mathematics, such as graph theory, group theory, stochastic processes, mathematical physics etc.
Since evolution algebras are not defined by identities, they can not belong to any well-known classes of non-associative algebras, as Lie, alternative and Jordan algebras.

The foundation of evolution algebra theory and applications in non-Mendelian genetics and Markov chains are
developed, with pointers to some further research topics was given in book \cite{TIAN2}.

In this paper, we study some properties of finite dimensional complex evolution algebras. Since any evolution algebra in a natural
basis is defined by a quadratic matrix, we study the connection between the algebraic structure of evolution algebras and matrices. More precise results are obtained for evolution algebras with non-singular matrices. For example, the only automorphisms for such algebras are the composition of basis permutation and the multiplication of basic vectors to scalars. Since in the matrix theory the Jordan form of the matrix is essential topic, we investigate a class of evolution algebras isomorphic to evolution algebras with Jordan form matrices. Thus we can distinguish the class of evolution algebras with a matrix in which the eigenvalues are known. Therefore, corresponding algebras can be investigated by the eigenvalues in algebraical point of view. Namely, the problem of reconstruction of Markov chains on trees \cite{MOS} which depends on the second eigenvalue can be studied by above evolution algebras.

In \cite{ETH3} it was pointed out for general genetic algebras that the nilpotent property is essential to these algebras and the definition as train algebras and baric algebras were formulated. By this means, we define nil, solvable, right-nilpotent and nilpotent evolution algebras as in \cite{CLOR} and study some properties of $n-$dimensional nilpotent evolution algebras. The notions as right nilpotency and nility for finite dimensional evolution algebras are equivalent \cite{CLOR}. In this work, we prove that any $n-$dimensional right-nilpotent evolution algebra is nilpotent. Moreover, for evolution algebras of dimension $n$ we describe some possible values for indexes of nilpotency and prove that $1+2^{n-1}$ is a maximal nilpotency index.

In \cite{TIAN2} the relation between graph theory and evolution algebras was given. The last section of this work is dedicated to the study of some evolution algebras defined by graphs, namely we find some algebraic properties of evolution algebras defined by complete and complete bipartite graphs and reformulate the graph planarity criteria in terms of evolution algebras.

\section{Preliminaries}
Now we define the main object of the paper.
\begin{defn}\cite{TIAN2} Let $E$ be a vector space over a field $K$ with defined multiplication
 $\cdot$ and a basis $\{e_1,e_2,\dots\}$ such that $$e_i\cdot e_j=0, \ i\neq j,$$
$$e_i\cdot e_i=\sum_k a_{ik}e_k, \ i \geq 1,$$ then $E$ is called evolution algebra and
basis $\{e_1,e_2,\dots\}$ is said to be natural basis. \end{defn}

From the above definition it follows that evolution algebras are commutative (therefore, flexible).

Let $E$ be a finite dimensional evolution algebra with natural
basis $\{e_1,\dots, e_n\},$ then
$$e_i\cdot e_i=\sum\limits_{j=1}^n a_{ij}e_j, \ 1\leq i \leq n,$$
where remaining products are equal to zero.

The matrix $A=(a_{ij})_{i,j=1}^{n}$ is called matrix of the algebra
$E$ in natural basis $\{e_1,\dots ,e_n\}.$

In \cite{TIAN2} conditions for basis transformations that preserve naturalness of the basis are given. Also, the relation between the matrices in a new and old natural basis is established in terms of new defined operation on matrices. Since this approach is not practical for our further purposes, below we give the following brief version in terms of its matrix elements.

Now let us consider non-singular linear transformation of the given natural basis $\{ e_1,\dots, e_n\}$ by matrix $T=(t_{ij})_{i,j=1}^n :$
$$f_i= \sum_{j=1}^n t_{ij}e_j, \ 1\leq i \leq n.$$

This transformation is isomorphism if and only if $f_i  f_j=0$
for all $i\neq j.$

Thus, $$\displaystyle f_i\cdot f_j=
\sum_{p=1}^n t_{ip}t_{jp} (e_p\cdot  e_p)=
\sum_{k=1}^n\left(\sum_{p=1}^n t_{ip}t_{jp} a_{pk}\right)e_k=0.$$

Hence, if $T$ is an isomorphism, then for $i\neq j$ and $1\leq k \leq n$ we have
$$\displaystyle\sum_{p=1}^n t_{ip}t_{jp} a_{pk}=0. \eqno(2.1)$$

Observe that
$$f_i\cdot f_i=
\sum_{p=1}^n t_{ip}^2 (e_p\cdot e_p)= \sum_{p=1}^n t_{ip}^2\sum_{k=1}^n
a_{pk}e_k=\sum_{k=1}^n\left(\sum_{p=1}^n t_{ip}^2
a_{pk}\right)e_k.$$

Now let $T_{ij}$ be the elements of matrix $T^{-1}.$ Then $\displaystyle e_k=\sum_{s=1}^n T_{ks}f_s$ and
$$f_i \cdot f_i=\sum_{k=1}^n\left(\sum_{p=1}^n t_{ip}^2
a_{pk}\right)\sum_{s=1}^nT_{ks}f_s=\sum_{s=1}^n
\left(\sum_{k=1}^n\sum_{p=1}^n t_{ip}^2 a_{pk}T_{ks}\right)f_s.$$

Hence, for the elements of the matrix $B=(b_{is})_{i,s=\overline{1,n}}$ of evolution algebra $E$ in natural basis $\{f_1,\dots,f_n\}$ we have
$$b_{is}= \sum_{k=1}^n\sum_{p=1}^n t_{ip}^2 a_{pk}T_{ks}. \eqno (2.2)$$

\begin{defn} An element $a$ of evolution algebra $E$ is called nil if there exists $n(a)\in \mathbb{N}$ such that $(\dots(\underbrace{(a\cdot a)\cdot a)\cdot \dots )\cdot a)}_{n(a)\textrm { times }}=0.$ Evolution algebra $E$ is called nil if any element of the algebra is nil.
\end{defn}

We introduce the following sequences:

 $$E^{(1)}=E,\, E^{(k+1)}=E^{(k)}E^{(k)},\, k \geq 1$$

$$E^{<1>}=E,\, E^{<k+1>}=E^{<k>}E,\, k \geq 1$$

$$ E^{1}=E,\,
E^{k}=\sum_{i=1}^{k-1}E^iE^{k-i},\, k \geq 1$$

Note that is not difficult to prove the following inclusions for $k\geq 1:$
$$E^{<k>}\subseteq E^k,\, E^{(k+1)}\subseteq E^{2^k}.$$

Also, note that since $E$ is commutative algebra we obtain
$\displaystyle E^{k}=\sum_{1\leq i\leq k-i}E^iE^{k-i}.$

\begin{defn} An evolution algebra $E$ is called

(i) solvable if there exists $n\in \mathbb{N}$ such that
$E^{(n)}=0$ and the minimal such number is called index of solvability;

(ii) right nilpotent if there exists $n\in \mathbb{N}$ such that $E^{<n>}=0$ and the minimal such number is called index of right nilpotency;

(iii) nilpotent if there exists $n\in \mathbb{N}$ such that
$E^{n}=0$ and the minimal such number is called index of nilpotency.
\end{defn}

Observe that if evolution algebra is nilpotent, then it is right nilpotent and solvable. The following example shows that solvable evolution algebra is not necessarily a right nilpotent algebra.

\begin{ex} Let $E$ be an evolution algebra with natural basis $\{e_1,\dots e_n\}$ and the following multiplication:
$$e_i e_i=e_1+\dots+e_n, \, 1\leq i \leq n-1$$
$$e_ne_n=(1-n)(e_1+\dots+e_n).$$
Then $E^{(3)}=0,$ but $E^{k}=\langle e_1+\dots +e_n\rangle$ for $k\geq 2.$
\end{ex}

The example described above in fact is a particular case of the following

\begin{pr}{pro-sol} Let $E$ be an $n-$dimensional complex  evolution algebra such that $\dim E^{(2)}=1.$ Then $E^{(3)}=0$ if and only if $E$ is isomorphic to an evolution algebra with natural basis $\{e_1,\dots,e_n\}$ with the following multiplication: $$e_ie_i=\lambda_i (e_1+\dots+e_k),\, 1\leq i \leq n,$$ where $\displaystyle\lambda_i \in \mathbb{C},\,\sum_{j=1}^k\lambda_j=0,\, \sum\limits_{j=1}^n|\lambda_j|^2\neq 0$ and $1\leq k \leq n.$
\end{pr}
\begin{proof} Since $E^{(2)}=1$ and $E^{(2)}$ is spanned by $e_ie_i,\, 1\leq i \leq n$ we obtain that they are collinear to a non-zero vector $x=a_1e_1+\dots+a_ne_n$. With the suitable natural basis change, one can assume that $x=e_1+\dots+e_k$ for some $1\leq k \leq n.$

Let $e_ie_i=\lambda_i x,\, 1\leq i \leq n$ and $\displaystyle \sum_{j=1}^n|\lambda_j|^2\neq 0.$ Then $E^{(3)}$ is spanned by $$xx=(e_1+\dots+e_k)^2=\sum_{j=1}^k \lambda_j x.$$

Hence, $E^{(3)}=0$ if and only if $\displaystyle \sum_{j=1}^k\lambda_j=0.$ \end{proof}

\begin{rem} Actually, the multiplication obtained in Proposition \ref{pro-sol} can be divided into two disjoint classes. First one, when $\lambda_i=0$ for all $1\leq i \leq k,$ then this evolution algebra is nilpotent. The second one is when $\lambda_i \neq 0$ for some $1\leq i \leq k.$ Then by natural basis transformation one can assume that $e_1e_1=e_1+\dots+e_k$ and hence, this evolution algebra is not nilpotent.
\end{rem}

In \cite{CLOR} the equivalence of right nilpotency and nility for finite dimensional complex evolution algebras is proved.

\begin{thm}\label{rnil} The following statements are equivalent:

a) The matrix of an evolution algebra $E$ can be transformed by natural basis permutation to

$$A=\left(\begin{array}{ccccc}
0 &a_{12}&a_{13}&\dots & a_{1n}\\
0 &0&a_{23}&\dots & a_{2n}\\
0 &0&0&\dots & a_{3n}\\
\vdots &\vdots&\vdots&\ddots & \vdots\\
0 &0&0&\dots & 0\\
\end{array}\right); \eqno (2.3)$$

b) Evolution algebra $E$ is right nilpotent algebra;

c) Evolution algebra $E$ is nil algebra.

\end{thm}

\section{Isomorphisms}

In case of evolution algebras with non-singular evolution matrices the problem of finding isomorphic algebras to the given one can be solved more precisely.

Let $E$ be an evolution algebra with matrix $A$ such that $\det A\neq 0.$
\begin{pr}\label{trans} $Aut(E)=\{ T_{\pi} \, | \, \pi \in S_n\},$  where
$T_{\pi}=(t_{ij})_{1\leq i,j\leq n}$ such that $t_{ij}\neq 0$ if and only if $j=\pi(i).$ Moreover, if $T_{\pi}$ is an automorphism of evolution algebra $E$ and $B=(b_{ij})_{1\leq i,j\leq n}$ is the matrix of $E$ in basis $T_{\pi}(e_1),\dots,T_{\pi}(e_n)$ then
$$b_{ij}=\frac{t_{i,\pi(i)}^2}{t_{j,\pi(j)}}\cdot a_{\pi(i)\pi(j)}. \eqno (3.1) $$
\end{pr}

\begin{proof}
Consider $(2.1)$ as a linear homogeneous system of equations in terms of unknowns $t_{i1}t_{j1},\dots, t_{in}t_{jn}.$
If $A$ is a non-singular matrix then from $(2.1)$ we obtain
$$\left\{%
\begin{array}{c}
    t_{i1}t_{j1}=0 \\
    t_{i2}t_{j2}=0 \\
    \dots  \\
    t_{in}t_{jn}=0 \\
\end{array}%
\right.$$ where $i\neq j.$

Since matrix $T$ is non-singular, in every row there is at least one non-zero element. But for any non-zero element $t_{ip}$ (in the $i-$th row) we have $t_{ip}t_{jp}=0$ for all $j\neq i.$
Therefore, $t_{jp}=0$ for $j\neq i.$ Now if for some $m\neq p$ we
have $t_{im}\neq 0,$ then similarly, we obtain $t_{jm}$ for all
$j\neq m.$ But this contradicts to non-singularity of matrix $T.$
Therefore, in every row and every column we have exactly one non-zero element, i.e., the matrix $T$ has the form described in the statement of the proposition.

Note that $\det T=(-1)^{\sigma(\pi)} t_{1\pi(1)}t_{2\pi(2)}\cdots t_{n\pi(n)},$ where $\sigma(\pi)$ is a signature of $\pi.$

We obtain that the group of automorphisms of $E$ is $\{ T_{\pi} \, | \, \pi \in S_n\}$ and  $T_{\pi} \circ T_{\tau}=T_{\tau \circ \pi}.$

Let us fix one $\pi \in S_n.$ Then $T(e_i)=t_{i,\pi(i)}e_{\pi(i)}$ for all $1\leq i \leq n.$ Now

$\displaystyle T(e_i) \cdot T(e_i)=t_{i,\pi(i)}^2 (e_{\pi(i)}\cdot e_{\pi(i)})=t_{i,\pi(i)}^2  \sum_{k=1}^n a_{\pi(i)k}e_k=$
$$t_{i,\pi(i)}^2  \sum_{k=1}^n a_{\pi(i)\pi(k)}e_{\pi(k)}=
\sum_{k=1}^n \frac{t_{i,\pi(i)}^2}{t_{k,\pi(k)}} a_{\pi(i)\pi(k)}T(e_k).$$

Hence, the elements of evolution matrix $B=(b_{ij})_{i,j=\overline{1,n}}$ of isomorphic algebra to $E$ satisfy $(3.1).$
\end{proof}

For a $\pi\in S_n$ denote by $s_{\pi}:\{1,2,\dots,n\} \setminus \{\pi^{-1}(n)\}\to\{1,2,\dots,n\}$ a one-to-one mapping defined by $s_{\pi}(i)=\pi^{-1}(1+\pi(i)).$

\begin{pr} Let $A=(a_{ij})_{1\leq i,j\leq n}$ be a matrix of an evolution algebra isomorphic to an evolution algebra with Jordan cell matrix with non-zero eigenvalue $\lambda.$ Then the only non-zero elements of $A$ are the diagonal elements and $a_{i,s_{\pi}(i)}$ for all
$i\neq \pi^{-1}(n)$ and
$\lambda=\displaystyle\frac{a_{ii}^2}{a_{i,s_{\pi}(i)}a_{s_{\pi}(i)s_{\pi}(i)}}$ for all $i\neq
\pi^{-1}(n).$
\end{pr}

\begin{proof}
First consider the isomorphism of evolution algebra with Jordan cell matrix with non-zero eigenvalue $\lambda.$ Since the matrix is non-singular, by the proof of Proposition \ref{trans} we obtain that it is in the form $T_{\pi}.$

For fixed $\pi \in S_n$ we put $T_{\pi}(e_i)=f_i$ and derive
$$f_i \cdot f_i=t_{i,\pi(i)}\lambda f_i+ \frac{t_{i,\pi(i)}^2}{t_{s_{\pi}(i),\pi(s_{\pi}(i))}}f_{s_{\pi}(i)}.$$

Hence the matrix of the new evolution algebra is a sum of non-singular diagonal matrix and a matrix that has exactly one non-zero element on each row except the $\pi^{-1}(n)-$th, which is a zero row and at most one non-zero element in each column.

Now let us fix a permutation $\pi \in S_n$ and consider matrix
$A=(a_{ij})_{i,j=1}^n$ with zero elements except the diagonal elements and $a_{i,s_{\pi}(i)}$ for all
$i\neq \pi^{-1}(n)$ and $s_{\pi}(i)=\pi^{-1}(1+\pi(i)).$ If this evolution algebra is isomorphic to an evolution algebra with Jordan cell matrix with eigenvalue $\lambda$ then $a_{ii}=\lambda t_{i, \pi(i)}$ for all $1\leq i \leq n$
and $a_{i,s_{\pi}(i)}=\frac{t_{i,\pi(i)}^2}{t_{s_{\pi}(i),\pi(s_{\pi}(i))}}.$

Since $t_{i,\pi(i)}= \frac1{\lambda}a_{ii}$ and
$t_{s_{\pi}(i),\pi(s_{\pi}(i))}=\frac1{\lambda}a_{s_{\pi}(i)s_{\pi}(i)}$ we obtain
$$a_{i,s_{\pi}(i)}=\frac{a_{ii}^2}{\lambda^2}\cdot
\frac{\lambda}{a_{s_{\pi}(i)s_{\pi}(i)}}=
\frac1\lambda\cdot\frac{a_{ii}^2}{a_{s_{\pi}(i)s_{\pi}(i)}} \textrm { and hence }\lambda= \frac{a_{ii}^2}{a_{i,s_{\pi}(i)}a_{s_{\pi}(i)s_{\pi}(i)}}.$$

Hence, if matrix $A$ satisfies
$\lambda=\frac{a_{ii}^2}{a_{i,s_{\pi}(i)}a_{s_{\pi}(i)s_{\pi}(i)}}$ for all $i\neq
\pi^{-1}(n),$ then evolution algebra with matrix $A$ is isomorphic to evolution algebra with Jordan cell matrix with eigenvalue $\lambda.$ This isomorphism has the matrix which is the inverse to $T=(t_{ij})_{i,j=1}^n,$ where $t_{i\pi(i)}=\frac1\lambda
a_{ii}$ and zero otherwise.
\end{proof}

The above result can be generalized to the case of Jordan form matrices. Let $J=J_1\oplus J_2 \oplus \dots \oplus J_r,$ where $J_i$ are Jordan cells of dimension $n_i$ with non-zero eigenvalue $\lambda_i.$

Now let us denote $\displaystyle \mu_k=\lambda_i$ for $n_1+\dots+n_{i-1}+1\leq k \leq n_1+\dots+n_i,\,1\leq i \leq r.$

Take $\pi\in S_n$ and denote  $s'_{\pi}:\{1,\dots,n\}\setminus \{\pi^{-1}(n_1),\dots,\pi^{-1}(n_r)\} \to \{1,\dots,n\}$  a one-to-one mapping defined by $s'_{\pi}(i)=\pi^{-1}(1+\pi(i)).$

\begin{cor}
Let $A=(a_{ij})_{1\leq i,j\leq n}$ be a matrix of an evolution algebra isomorphic to an evolution algebra with Jordan form matrix $J.$ Then the only non-zero elements of $A$ are the diagonal elements and $a_{i,s'_{\pi}(i)}$ such that
$\displaystyle\frac{a_{ii}^2}{a_{i,s'_{\pi}(i)}a_{s'_{\pi}(i)s'_{\pi}(i)}}=\frac{\mu_i^2}{\mu_{s'_{\pi}(i)}}$  for \\ $i\not \in \{\pi^{-1}(n_1),\dots,\pi^{-1}(n_r)\}.$
\end{cor}

\section{Nilpotency of evolution algebras}

Let us now consider an evolution algebra $E$ with Jordan cell with eigenvalue $\lambda$.

\begin{pr} If $\lambda \neq 0$ then $E$ is neither
solvable nor right nilpotent and therefore is not nilpotent.\end{pr}

\begin{proof} Since $\lambda \neq 0$ then evolution matrix is
non-degenerated. Therefore, $E^2=E^{(2)}=E^{<2>}=E.$ By simple
induction we obtain $E^k=E^{(k)}=E^{<k>}=E$ and the statement of
the proposition is verified. \end{proof}

\begin{pr} \label{jordan-nil} For an evolution algebra with Jordan cell matrix and eigenvalue $\lambda=0$ the following statements hold:

(i) $E$ is one generated;

(ii) $E$ is solvable with index of solvability $n+1;$

(iii) $E$ is right nilpotent with index of right nilpotency $n+1;$

(iv) $E$ is nilpotent with index of nilpotency $2^{n-1}+1.$
\end{pr}

\begin{proof} (i) From $\lambda=0$ it follows that for basis elements
$e_i$ of $E$ we have $[e_i,e_i]=e_{i+1}$ for all $1\leq i \leq n-1
$ and $[e_n,e_n]=0.$

Therefore, $E$ is one-generated: $E=id\langle e_1\rangle.$\\

(ii) First observe that $E^{(2)}=\langle e_2,\dots,e_n\rangle.$

If for some $k$ we have
$$E^{(k)}=\langle e_k,e_{k+1},\dots,e_n \rangle,$$
then for $k+1$ we obtain

$$E^{(k+1)}=E^{(k)}E^{(k)}=\langle e_{k+1},\dots,e_n \rangle. $$

Therefore, $E^{(n)}=\langle e_n \rangle$ and $E^{(n+1)}= 0$ and
(ii) is  verified.\\

(iii) is similar to (ii).\\

(iv) We claim that
$$E^{2^k+1}=E^{2^k+2}=\dots=E^{2^{k+1}}=\langle
e_{k+2},\dots,e_n\rangle $$ for all $0\leq k \leq n-2.$

Indeed, for $k=0$ we have $E^2=EE=\langle e_2,\dots,e_n\rangle.$

For $k=1$ we have $$E^{2+1}=E^3=EE^2=\langle
e_3,\dots,e_n\rangle,$$
$$E^{2^2}=E^4=EE^3+E^2E^2=\langle e_3,\dots,e_n\rangle.$$

Assume that
$$E^{2^{k-1}+1}=E^{2^{k-1}+2}=\dots=E^{2^k}=\langle
e_{k+1},\dots,e_n\rangle. $$

Using this assumption we obtain
$$E^{2^k+1}=EE^{2^k}+E^2E^{2^k-1}+\dots+E^{2^{k-1}}E^{2^{k-1}+1}=$$
$$EE^{2^k}+E^2E^{2^k}+\dots+E^{2^{k-1}}E^{2^k}=$$
$$(E+E^2+E^3+\dots+E^{2^{k-1}}) E^{2^k} =EE^{2^{k}}=\langle e_{k+2},\dots,e_n\rangle. $$

Also
$$E^{2^{k+1}}=EE^{2^{k+1}-1}+E^2E^{2^{k+1}-2}+\dots+E^{2^{k}}E^{2^k}\supseteq$$
$$E^{2^{k}}E^{2^k}=\langle e_{k+2},\dots,e_n\rangle. $$
So we obtain
$$\langle e_{k+2},\dots,e_n\rangle=E^{2^k+1}\supseteq
E^{2^k+2}\supseteq\dots\supseteq E^{2^{k+1}}\supseteq
\langle e_{k+2},\dots,e_n\rangle. $$

Hence, $$E^{2^k+1}= E^{2^k+2}=\dots=E^{2^{k+1}}=\langle
e_{k+2},\dots,e_n\rangle.$$

Therefore, $E^{2^{n-1}}=\langle e_n \rangle $ and
$E^{2^{n-1}+1}=0.$

Hence, $E$ is nilpotent with nilpotency index equal to
$1+2^{n-1}.$ \end{proof}

\begin{rem}We should note that the statements $(ii)-(iv)$ of Proposition \ref{jordan-nil} are equivalent, since one can show that each of them is equivalent to $\lambda=0.$ However, statement $(i)$ is not equivalent to $\lambda =0$ since for $\lambda =1$ one can prove that $E$ is generated by the element $e_1+e_2.$
\end{rem}

Observe that any evolution subalgebra of an evolution algebra is an ideal.
Therefore if we consider an evolution algebra $E_J$ with matrix
$J$ in Jordan form $J= J_1 \oplus J_2\oplus \dots \oplus J_r $ where $J_i$ are Jordan cells of dimension $n_i$ with eigenvalues $\lambda_i,$ then   $$E_J=E_1\oplus E_2 \oplus \dots \oplus E_r,$$ where $E_i=\langle e_{n_{i-1}+1},\dots,e_{n_i}\rangle.$

Now we have $E_J^k=E_1^k\oplus E_2^k\oplus \dots \oplus E_r^k$ and
therefore $E_J$ is nilpotent if and only if every $E_i$ is
nilpotent. Since we have obtained the criteria of nilpotency of
Jordan blocks, we obtain

\begin{cor}\label{cor-jordan} $E_J$ is nilpotent (with index of nilpotency equal to $\max_{1\leq i \leq r} \{1+2^{n_i-1}\}$) if and only if $J$ has only zero eigenvalues. The same assertion holds for right nilpotency and solvability with corresponding indexes equal to $1+\max_{1\leq i \leq r}\{n_i\}.$
\end{cor}

Note that from the Corollary \ref{cor-jordan} it follows that for every $1\leq k\leq n$ we obtain an example of nilpotent evolution algebra with index of nilpotency equal to $1+2^{k-1}.$

The following theorem represents the criteria of nilpotency of finite dimensional evolution algebra.

\begin{thm} \label{nil-index} Let $E$ be an $n-$dimensional evolution algebra. Then $E$ is nilpotent if and only if the matrix of evolution algebra $A$ can be transformed by the natural basis permutation to form (2.3). Moreover, the index of nilpotency of evolution algebra $E$ is not greater then $2^{n-1}+1.$
\end{thm}

\begin{proof} Let $E$ be a nilpotent. Then it is right nilpotent and therefore, by Theorem \ref{rnil} the matrix of this evolution algebra can be transformed by the natural basis permutation to from $(2.3).$

Now let the matrix $A$ of $E$ can be transformed by the natural basis permutation to form $(2.3).$

Assume that $a_{12}a_{23}\dots a_{n-1n} \neq 0.$ Similar to the proof of $(iv)$ in Proposition \ref{jordan-nil} one can verify
$$E^{2^k+1}=E^{2^k+2}=\dots=E^{2^{k+1}}=\langle
e_{k+2},\dots,e_n\rangle $$ for all $0\leq k \leq n-2.$

Therefore, $E^{2^{n-1}}=\langle e_n \rangle $ and
$E^{2^{n-1}+1}=0.$

Hence, $E$ is nilpotent with nilpotency index equal to
$1+2^{n-1}.$\\

Now assume that $a_{12}a_{23}\dots a_{n-1n}= 0.$ In this case we claim that
$$\langle
e_{k+2},\dots,e_n\rangle \supseteq E^{2^k+1} $$ for all $0\leq k \leq
n-2.$

Indeed, for $k=0$ we have $E^2=EE \subseteq \langle
e_2,\dots,e_n\rangle.$

For $k=1$ we have $$E^{2+1}=E^3=EE^2\subseteq\langle
e_3,\dots,e_n\rangle.$$

Assume that
$$\langle
e_{k+1},\dots,e_n\rangle\supseteq E^{2^{k-1}+1}. $$

Using this assumption we obtain
$$E^{2^k+1}=EE^{2^k}+E^2E^{2^k-1}+\dots+E^{2^{k-1}}E^{2^{k-1}+1} \subseteq $$
$$EE^{2^{k-1}+1}+E^2E^{2^{k-1}+1}+\dots+E^{2^{k-1}}E^{2^{k-1}+1}=$$
$$(E+E^2+E^3+\dots+E^{2^{k-1}}) E^{2^{k-1}+1} \subseteq EE^{2^{k-1}+1} \subseteq \langle e_{k+2},\dots,e_n\rangle. $$

So we obtain
$$\langle e_{k+2},\dots,e_n\rangle \supseteq E^{2^k+1}.$$

Therefore, $\langle e_n \rangle \supseteq E^{2^{n-2}+1}.$

Hence,
$$E^{2^{n-1}+1}=EE^{2^{n-1}} + E^2 E^{2^{n-1}-1} \dots + E^{2^{n-2}}E^{2^{n-2}+1} \subseteq$$
$$(E+E^2+\dots E^{2^{n-2}}) \langle e_n \rangle \subseteq E \langle e_n \rangle =0.$$

Thus, $E$ is nilpotent with nilpotency index not greater then
$1+2^{n-1}.$
\end{proof}

\begin{cor} For finite dimensional complex evolution algebra notions as nil, nilpotent and right nilpotent algebras are equivalent. However, the indexes of nility, right nilpotency and nilpotency do not coincide in general.
\end{cor}

The following proposition excludes significantly many possible values that a nilpotency indexes of $n-$dimensional evolution algebras can take.

\begin{pr} \label{index-great} Let $E$ be a nilpotent evolution algebra with index of nilpotency not equal to  $2^{n-1}+1.$ Then it is not greater then $2^{n-2}+1.$
\end{pr}
\begin{proof} Since $E$ is nilpotent, we assume that the matrix $A$ of $E$ in the natural basis $\{e_1,\dots,e_n\}$ is in the form $(2.3).$ 

From the proof of Proposition \ref{nil-index} it follows that
$a_{12}a_{23}\dots a_{n-1n}= 0$ and
$$\langle
e_{k+2},\dots,e_n\rangle \supseteq E^{2^k+1}$$ for all $0\leq k \leq
n-2.$

Assume that $E$ is nilpotent with index of nilpotency greater then
$2^{n-2}+1$ and not equal to $2^{n-1}+1.$

Then $\langle e_n \rangle \supseteq E^{2^{n-2}+1}$ and since
$E^{2^{n-2}+1}\neq 0$ we obtain $E^{2^{n-2}+1}=\langle
e_n \rangle.$

Therefore, $\langle e_{n-1},e_n \rangle \supseteq
E^{2^{n-3}+1}\supseteq E^{2^{n-3}+2}\supseteq \dots \supseteq
E^{2^{n-2}}\supseteq \langle e_n  \rangle.$

Now if $E^{2^{n-3}+1}= E^{2^{n-3}+2}= \dots=E^{2^{n-2}}= \langle
e_n  \rangle$ then
$$E^{2^{n-2}+1}=EE^{2^{n-2}} + E^2 E^{2^{n-2}-1} \dots +
E^{2^{n-3}}E^{2^{n-3}+1} \subseteq$$
$$(E+E^2+\dots +E^{2^{n-3}}) \langle e_n \rangle = E \langle e_n \rangle
=0$$ which is a contradiction. Hence, $\langle e_{n-1},e_n \rangle
= E^{2^{n-3}+1}.$

Now assume that $\langle e_{n-k},\dots, e_n \rangle =
E^{2^{n-k-2}+1}.$

Then $$\langle e_{n-k-1},e_{n-k},\dots, e_n \rangle \supseteq
E^{2^{n-k-3}+1}\supseteq E^{2^{n-k-3}+2}\supseteq \dots \supseteq
E^{2^{n-k-2}}\supseteq \langle e_{n-k},\dots, e_n  \rangle.$$ If
$E^{2^{n-k-3}+1} \neq \langle e_{n-k-1},e_{n-k},\dots, e_n
\rangle$ then $E^{2^{n-k-3}+1}=E^{2^{n-k-3}+2}=\dots =
E^{2^{n-k-2}}=\langle e_{n-k},\dots, e_n \rangle$ and
$$E^{2^{n-k-2}+1}=EE^{2^{n-k-2}}+\dots+E^{2^{n-k-3}}E^{2^{n-k-3}+1}=E \langle e_{n-k},e_{n-k},\dots, e_n
\rangle \subseteq \langle e_{n-k+1},e_{n-k},\dots, e_n \rangle$$
which contradicts to $\langle e_{n-k},\dots, e_n \rangle =
E^{2^{n-k-2}+1}.$

Hence this assumption is true and therefore $E^2 =\langle
e_2,\dots e_n \rangle $ which is also a contradiction to
$a_{12}\dots a_{n-1 n}=0.$
\end{proof}

The following example shows that there exist evolution algebras with index of nilpotency greater then $1+2^{k-3}$ and less then $1+2^{k-2}$ for all $4\leq k \leq n.$

\begin{ex} Consider an evolution algebra $E_k$ with basis $\{e_1,\dots, e_n\}$ and the following multiplication table:
$$\begin{array}{ll}
    e_1 e_1=e_2+e_3+\dots+e_k,&\\
    e_2 e_2=-e_4,&\\
    e_i e_i =e_{i+1},\,\, 3\leq i\leq k-1 & \textrm{and } 4\leq k\leq n.\\
\end{array}$$
Then one can show that $E_k^{3 \cdot 2^i}=\langle e_{4+i},\dots,e_k\rangle$ for $0\leq i\leq k-4$ and index of nilpotency of this algebra is $1+3\cdot 2^{k-4}.$
\end{ex}

Now we will consider a nilpotent evolution algebra with matrix $(2.3)$ and a condition $\dim E^2=n-2.$ Then $rank A=n-2.$ This implies that there are $1\leq i\leq n-1$ and $2\leq j\leq n$ such that $i-$th row is linear dependent to other rows and $j-$th column is linear dependent to other columns.

\begin{pr} Let $\dim E^2=n-2$ and $i-$th row $(1\leq i\leq n-1)$ is linear dependent to other rows and $j-$th column $(2\leq j\leq n)$ is linear dependent to other columns. Then
$$\dim E^3=
\left \{
\begin{array}{ccl}
n-3& \textrm{ if }& i=1 \textrm{ or }\\
&& j=n \textrm{ or }\\
&& j\neq n \textrm{ and } j-\textrm{th column is non-zero or }\\
&& i\neq 1, j=i \textrm{ and } j-\textrm{th column is zero } \\
n-4& \textrm{ if }& i\neq 1,j\neq i,n \textrm{ and } j-\textrm{th column is zero }\\
\end{array}\right.$$
Moreover, $i=1$ implies $j=2$ and $j=n$ implies $i=n-1.$
\end{pr}

\begin{proof}
Consider $$E^3=E\cdot E^2=\left\langle
\begin{array}{rrrc}
    e_2(e_1 e_1),&e_3(e_1 e_1),&\dots &e_{n-1}(e_1 e_1),\\
    &e_3(e_2 e_2),&\dots&e_{n-1}(e_2 e_2),\\
            & &\ddots &\vdots\\
    &&&e_{n-1}(e_{n-2} e_{n-2})\\
    \end{array} \right\rangle=$$

$$\left\langle
\begin{array}{rrrc}
  a_{12}(e_2 e_2),&a_{13}(e_3 e_3),&\dots& a_{1\,n-1}(e_{n-1} e_{n-1}), \\
    &a_{23}(e_3 e_3),&\dots & a_{2\,n-1}(e_{n-1} e_{n-1}),\\
        & &\ddots &\vdots\\
    &&&a_{n-2\, n-1}(e_{n-1} e_{n-1})\\
    \end{array} \right\rangle \eqno (4.1)$$

    Denote by $L:= \langle e_2e_2,  e_3e_3, \dots,  e_{n-1}e_{n-1}\rangle.$ Obviously, $E^3\subseteq L.$

    If $i=1$ then we obtain $a_{12}=0$ and $a_{23}\dots a_{n-1 n}\neq 0.$ Hence, $E^2=\langle e_3,\dots, e_n \rangle$ and $E^3=\langle e_4,\dots, e_n \rangle.$ Moreover, applying the same arguments as in the proof of Proposition \ref{jordan-nil} $(i)$ we obtain
    $$E^{2^{k-1}+1}=E^{2^k+2}=\dots=E^{2^{k}}=\langle
e_{k+2},\dots,e_n\rangle $$ for all $1\leq k \leq n-2$ and index of nilpotency for this algebra in this case is $1+2^{n-2}.$

Now let $2\leq i \leq n-1.$ Then $\dim L=n-3.$

If $j=n,$ then $a_{n-1 n}=0$ and therefore, $i=n-1.$ Hence, from $(4.1)$ one obtains $E^3=\langle e_2e_2,   e_3e_3, \dots,  e_{n-2}e_{n-2}\rangle.$ Thus, $\dim E^3=n-3.$

Now if $j\neq n$ and $j-$th column is non-zero column then one can easily see from $(4.1)$ that again $E^3=L.$ Hence, $\dim E^3=n-3.$

If  $j\neq n$ and $j-$th column is zero column, then $$E^3=\langle(e_2e_2),\dots, (e_{j-1}e_{j-1}), (e_{j+1}e_{j+1}),\dots, (e_{n-1} e_{n-1})\rangle.$$

Now if $i=j$ then $E^3=L.$ In this case $\dim E^3=\dim L=n-3.$

If $i\neq j$ then $\dim E^3=\dim L-1.$

Hence, the statement of the proposition is verified.
\end{proof}

\section{Graphs and Evolution Algebras}
In this section we will try to transfer some properties of graphs to algebraic properties of evolution algebras defined by them. For definition of graphs and their properties see \cite{GRAPH}.

The definition of evolution algebra defined by graph and the next theorem was given in \cite{TIAN2} for simple graphs. However, one can easily formulate the analogous definition and prove the theorem for directed graphs.

\begin{defn} Let $D=(V,E)$ be a directed graph with $n$ vertices from the set $V,$ the sorted edges from the set $E$ and the adjacency matrix $A=(a_{ij})_{1\leq i,j\leq n}.$ Then evolution algebra determined by this graph is an algebra $E(D)=\langle e_1,\dots, e_n\rangle$ with the following multiplication: $$e_i\cdot e_j =0, \,\,\,\,\,\,\,\,\, e_i\cdot e_i =\sum_{k=1}^n a_{ik}e_k \textrm{ for all } 1\leq i \neq j\leq n.$$
\end{defn}

\begin{thm} If graphs $G_1$ and $G_2$ are isomorphic as graphs, then $E(G_1)$ and $E(G_2)$ are isomorphic as evolution algebras.
\end{thm}

For definitions of complete and complete bipartite graph and their properties see \cite{GRAPH}.

\begin{defn} Evolution algebra determined by a complete (complete bipartite) graph is called a compete (complete bipartite) evolution algebra.
\end{defn}

\begin{pr} Let $E$ be an evolution algebra with natural basis $\{e_1,\dots,e_n\}$ and with a matrix of non-negative integers. Then this algebra is a complete evolution algebra if and only if there exists $k \,(1\leq k \leq n)$ such that $$(\dots ((e_{i_1}  e_{i_1})e_{i_2})\dots ) e_{i_k}= e_1+\dots+e_{i_k-1}+e_{i_k+1}+\dots+e_n\eqno(5.1)$$ for $i_1,\dots , i_k \in \{1,\dots,n\}.$
\end{pr}

\begin{proof} Let $ \displaystyle e_i\cdot e_i= \sum_{k=1}^n a_{ik}e_k.$ Then from the condition of proposition we obtain $a_{i_1 i_2}a_{i_2 i_3}\dots a_{i_{k-1} i_k} (e_{i_k} e_{i_k}) =  e_1+\dots+e_{i_k-1}+e_{i_k+1}+\dots+e_n$ for all $i_1,\dots , i_k \in \{1,\dots,n\}.$

The last one implies $a_{i_1 i_2}a_{i_2 i_3}\dots a_{i_{k-1} i_k} a_{i_k p} =1$ for all $p\neq i_k$ and
$a_{i_1 i_2}a_{i_2 i_3}\dots a_{i_{k-1} i_k} a_{i_k i_k} =0.$

Therefore, $a_{i_k i_k}=0$ and since the elements of the matrix $A$ are non-negative integers we obtain $a_{i_k p}=1$ for $p\neq i_k.$
Since $i_k$ can take arbitrary values from $\{1,\dots,n\}$ we obtain that $a_{ij}=1$ for all $i\neq j$ and $a_{ii}=0.$ Thus this algebra is a complete evolution algebra.

The proof in the opposite direction is obvious.
\end{proof}

Let us denote $Z_1=\{1,\dots,n\},\,Z_2=\{n+1,\dots,2n\}$ and for a natural $q$ by $Z_q$ we mean $Z_1$ if $q$ is odd and $Z_2$ otherwise.

\begin{pr} Let $E$ be an evolution algebra with natural basis $\{e_1,\dots,e_{2n}\}$ and with a matrix of non-negative real elements. Then this algebra has the following matrix in the natural basis
$$\left(
\begin{array}{cccccc}
0&\dots&0&\alpha&\dots&\alpha\\
\vdots&&\vdots&\vdots&&\vdots\\
0&\dots&0&\alpha&\dots&\alpha\\
\beta&\dots&\beta&0&\dots&0\\
\vdots&&\vdots&\vdots&&\vdots\\
\beta&\dots&\beta&0&\dots&0\\
\end{array}\right),  \alpha\beta=1 \textrm{ if and only if }$$

$$(\dots ((e_{i_1}  e_{i_1})e_{i_2})\dots ) e_{i_k}=\left\{
\begin{array}{cl}
0,&\textrm{ there exists } p \textrm{ such that }\\
    &  \{e_{i_p},e_{i_{p+1}}\}\subseteq Z_1 \textrm { or } \{e_{i_p},e_{i_{p+1}}\}\subseteq Z_2\\
    \displaystyle\sum_{j\in Z_{i_k+1}} e_j,& \textrm{  otherwise  } \\
\end{array}\right. \eqno(5.2)$$
for $i_1,\dots , i_k,p \in \{1,\dots,2n\}$ and for some $1\leq k \leq 2n.$
\end{pr}

\begin{proof}  Let the condition of the proposition be true. Assume that $e_i\cdot e_i= \sum_{k=1}^{2n} a_{ik}e_k.$

Taking $i_1\in Z_q,i_2\in Z_{q+1},\dots,i_k\in Z_{q+k-1}$ for some natural $q$ one obtains
$a_{i_1i_2}\dots a_{i_{k-1}i_k}(e_{i_k}e_{i_k})\neq 0.$

Now taking $i_1,i_2\in Z_q,i_3\in Z_{q+1},$ $\dots,i_k\in Z_{q+k-2}$ for some natural $q$ one obtains
$a_{i_1i_2}\dots a_{i_{k-1}i_k}(e_{i_k}e_{i_k})=0.$

Hence, if $i_1,i_2 \in Z_q$ for some natural $q$ then $a_{i_1 i_2}=0.$

Also by taking $i_1\in Z_q,i_2\in Z_{q+1},\dots,i_k\in Z_{q+k-1}$ for some $q,$ we obtain
$a_{i_1i_2}\dots a_{i_{k-1}i_k}a_{i_kj}=1$ for all $j\in Z_{i_k+1}.$

Hence, $a_{i_k i}=0$ for all $i\in Z_{i_k}$ and $a_{i_k p}=\frac1{a_{i_1 i_2}a_{i_2 i_3}\dots a_{i_{k-1} i_k}}$ for all $p\in Z_{i_k+1}.$
Since $i_k$ can take arbitrary values from $Z_{i_k}$ we obtain that in each row the non-zero values of the elements are similar.

Now we can assume that $e_pe_p=c_p(e_{n+1}+\dots+e_{2n}),\,1\leq p \leq n$ and  $e_qe_q=c_q(e_{1}+\dots+e_{n}),\,n+1\leq q \leq 2n$ for some $c_1,\dots,c_{2n}\in \mathbb{C}.$

Taking $i_1\in Z_q,i_2\in Z_{q+1},\dots,i_k\in Z_{q+k-1}$ for some $q$ one obtains
$c_{i_k}=\frac1{c_{i_1}\dots c_{i_{k-1}}}.$ Since we can put every $i\in Z_{q+k-1}$ instead of $i_k$ we obtain that for $i\in Z_{q+k-1}$ we have $c_i=c_{i_k}.$

Analogously, if $i_1\in Z_{q+1},i_2\in Z_{q+2},\dots,i_k\in Z_{q+k}$ for some $q$ one obtains
$c_{i_k}=\frac1{c_{i_1}\dots c_{i_{k-1}}}.$ Since we can put every $i\in Z_{q+k}$ instead of $i_k$ we obtain that for $i\in Z_{q+k-1}$ we have $c_i=c_{i_k}.$ So the matrix of $A$ is as follows:

$$\left(
\begin{array}{cccccc}
0&\dots&0&\alpha&\dots&\alpha\\
\vdots&&\vdots&\vdots&&\vdots\\
0&\dots&0&\alpha&\dots&\alpha\\
\beta&\dots&\beta&0&\dots&0\\
\vdots&&\vdots&\vdots&&\vdots\\
\beta&\dots&\beta&0&\dots&0\\
\end{array}\right)$$

Now if $k$ is even then from $c_{i_1}\dots c_{i_{k-1}}c_{i_k}=1,$ where $c_{i_{2m-1}}\in Z_q$ and $c_{2m}\in Z_{q+1}$ for some $q$ we obtain $(\alpha\beta)^{\frac k2}=1.$ Hence, in this case $\alpha \beta=1.$

Now if $k$ is odd then from $c_{i_1}\dots c_{i_{k-1}}c_{i_k}=1,$ where $c_{i_{2m-1}}\in Z_q$ and $c_{2m}\in Z_{q+1}$ for some $q$ we obtain $(\alpha\beta)^{\frac {k-1}2} \alpha=1$ and for another set of $i_1,\dots,i_k$ we obtain $(\alpha\beta)^{\frac {k-1}2} \beta=1.$

This implies $\alpha=\beta=1.$ The proof in the opposite direction is obvious.
\end{proof}

\begin{cor}Let $E$ be an evolution algebra with natural basis $\{e_1,\dots,e_{2n}\}$ and with a matrix of non-negative integer elements.
Then this algebra is complete bipartite evolution algebra with partitions of equal size if and only if it satisfies (5.2)
with $i_1,\dots , i_k,p \in \{1,\dots,2n\}$ for some even $2\leq k \leq 2n.$
\end{cor}

Now we will define the concept in evolution algebra defined by graph which corresponds to a subgraph. In fact, by the renumbering we can always suppose that the vertex of a subgraph $G_1$ are $e_1,\dots e_k$ of the graph $G$ with vertices $e_1,\dots,e_k,\dots , e_n.$
The matrix corresponding to subgraph is a submatrix obtained by intersection of first $k$ rows and columns of $A.$ Finally, we obtain new evolution algebra in the basis $\{e_1,\dots e_k\}$ with corresponding matrix which is a submatrix of $A$ of size $k.$  Such type of evolution algebras we will denote by $E(G_1).$  In case when $\langle e_{k+1},\dots,e_n \rangle$ form an evolution subalgebra of $E$, then $E(G_1)$ is a quotient algebra of $E$ by $\langle e_{k+1},...,e_n \rangle.$

In graph theory, a planar graph is a graph that can be embedded in the plane, i.e., it can be drawn on the plane in such a way that its edges intersect only at their endpoints. Now we define planar evolution algebras.

\begin{defn} Evolution algebra determined by a planar graph is called a planar evolution algebra.
\end{defn}

In graph theory the process of the shrinkage of a graph plays an important role in the theory of planar graphs.

Let $G$ be a graph with vertices $e_1,\dots,e_n$ and adjacency matrix $A=(a_{ij})_{1\leq i,j\leq n}.$

If the vertices $e_p$ and $e_q$ are neighboring ($a_{pq}=a_{qp}\neq0$), then we can shrinkage $e_q$ to $e_p$ by the following way:

consider the graph with verteces  $e_1,\dots,e_p,\dots,e_{q-1},e_{q+1},\dots,e_n$ with an adjacency matrix $\widetilde{A}=(\widetilde{a_{ij}})_{1\leq i,j\leq n-1}$ which is obtained from the matrix $A$ with the following procedure:

first, we replace the elements $a_{pk}$ of the $p-$th row  by the elements $\max(a_{pk},a_{qk})$ and then eliminate the $q-$th row and the same column.

Evidently, the procedure is commutative under considering vertexes $e_p$ and $e_q.$ Therefore, we can always assume that $p<q.$

By the matrix $\widetilde{A}$ and basis $\widetilde{e_i}=e_i, 1\leq i \leq q-1$ and $\widetilde{e_j}=e_{j+1}, q\leq j \leq n-1$ we determine the evolution algebra $\widetilde{E}$ with natural basis $\{\widetilde{e_1},\dots,\widetilde{e}_{n-1}\}$ and table of multiplication:
$$\widetilde{e}_p^2=\sum_{k=1}^{q-1} \max(a_{pk},a_{qk})\widetilde{e}_k+
\sum_{k=q+1}^{n}\max(a_{pk},a_{qk})\widetilde{e}_{k-1},$$
$$\widetilde{e}_i^2=\sum_{k=1}^{p-1} a_{ik}\widetilde{e}_k+
\max(a_{pp},a_{pq})\widetilde{e}_p+
\sum_{k=p+1}^{q-1} a_{ik}\widetilde{e}_k+ \sum_{k=q+1}^{n} a_{ik}\widetilde{e}_{k-1}, \, 1\leq i\neq p \leq n-1.$$

In graph theory there are several criterias of planarity graphs. Now Harary-Tatta criteria states that a finite graph is planar if and only if it does not contain a subgraph that is a shrinkage of $K_5$ (the complete graph on five vertices) or $K_{3,3}$ (complete bipartite graph on six vertices, three of which connect to each of the other three). Due to this criteria and previous propositions we can reformulate this criteria in terms of evolution algebras defined by graphs.

\begin{thm} Finite dimensional evolution algebra is planar if and and only if there is no evolution algebra $E(G_1)$ for all $G_1$ subgraph of $G$ such that $E(G_1)$ can be shrinkage to 5 or 6-dimensional evolution algebras satisfying (5.1) or (5.2), respectively.
\end{thm}

\textbf{Acknowledgments.} This work is supported in part by the PAICYT, FQM143 of Junta de Andaluc\'ia (Spain). The third named author was supported by the grant NATO-Reintegration ref. CBP.EAP.RIG. 983169. The last named author would like to acknowledge ICTP OEA-AC-84 for a given support.

\end{document}